\documentclass[letterpaper, 10 pt, journal, twoside]{IEEEtran}

\IEEEoverridecommandlockouts

\usepackage{graphics,graphicx} % for pdf, bitmapped graphics files
\usepackage{epsfig} % for postscript graphics files
\usepackage{mathptmx} % assumes new font selection scheme installed
\usepackage{times} % assumes new font selection scheme installed
\usepackage{amsmath} % assumes amsmath package installed
\usepackage{amssymb,nth, mathtools,dsfont}  % assumes amsmath package installed
\usepackage{amsbsy}
\usepackage{epstopdf}
\usepackage{cite}
\usepackage{romannum}
\usepackage{xcolor}
\usepackage{subcaption}
\usepackage{multirow}

\newtheorem{theorem}{Theorem}
\newtheorem{lemma}{Lemma}
\newtheorem{definition}{Definition}
\newtheorem{remark}{Remark}
\newtheorem{conjecture}{Conjecture}
\newtheorem{corollary}{Corollary}

\makeatletter
\newcommand*{\rom}[1]{\expandafter\@slowromancap\romannumeral #1@}

\title{Noncausal FIR Zames-Falb Multiplier Search for Exponential Convergence Rate}

\author{Jingfan Zhang, Peter Seiler and Joaquin Carrasco
\thanks{Jingfan Zhang and Joaquin Carrasco are with the Control Systems Centre, %\newline
 School of Electrical and Electronic Engineering, University of Manchester, M13 9PL, UK. 
         {\tt\small jingfan.zhang@manchester.ac.uk}; 
         {\tt\small joaquin.carrascogomez@manchester.ac.uk}      
       }
   \thanks{Peter Seiler is with the Aerospace Engineering and Mechanics Department, University of Minnesota, Minneapolis, MN 55455 USA. 
   {\tt\small  seile017@umn.edu}   
}
}

\begin{document}

\maketitle
\thispagestyle{empty}
\pagestyle{empty}

\begin{abstract}
	In the existing literature, there are two approaches to estimate tighter bounds of the exponential convergence rate of stable Lur'e systems. On one hand, the classical integral quadratic constraint (IQC) framework can be applied under loop-transformation, so the stability of the new loop implies the convergence of the original loop. On the other hand, it is possible to modify the IQC framework, the so-called $\rho$-IQC framework, in such a way that the convergence rate is directly obtained over the original loop. 
	In this technical note, we extend the literature results from the search for a causal finite impulse response (FIR) Zames-Falb multiplier to the noncausal case. We show that the multipliers by the two approaches are equivalent by a change of variable. However,  the factorisation of the Zames-Falb $\rho$-IQC is restricted compared to the  Zames-Falb IQC, so a unified factorisation is proposed. 
	Finally, numerical examples illustrate that noncausal multipliers lead to less-conservative results.

\end{abstract}


\section{Introduction}

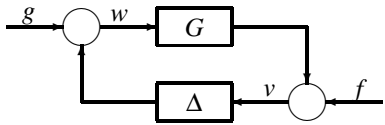
\begin{figure}[ht]
	\centering
	\ifx\JPicScale\undefined\def\JPicScale{1}\fi
\unitlength \JPicScale mm
\begin{picture}(50,15)(0,0)

\linethickness{0.3mm}
\put(10,2.5){\line(1,0){10}}
\put(10,2.5){\line(0,1){7.5}}
\put(10,10){\vector(0,1){0.12}}
\put(10,12.5){\circle{5}}

%\linethickness{0.5mm}
%\put(12,15){\line(1,0){2}}
\linethickness{0.3mm}
\put(0,12.5){\line(1,0){7.5}}
\put(7.5,12.5){\vector(1,0){0.12}}
\put(12.5,12.5){\line(1,0){7.5}}
\put(20,12.5){\vector(1,0){0.12}}
\put(30,12.5){\line(1,0){10}}
\put(40,5){\line(0,1){7.5}}
\put(40,5){\vector(0,-1){0.12}}

\put(40,2.5){\circle{5}}
\put(42.5,2.5){\line(1,0){7.5}}
\put(42.5,2.5){\vector(-1,0){0.12}}
\put(30,2.5){\line(1,0){7.5}}
\put(30,2.5){\vector(-1,0){0.12}}

\put(20,10){\line(1,0){10}}
\put(20,15){\line(1,0){10}}
\put(20,10){\line(0,1){5}}
\put(30,10){\line(0,1){5}}

\put(20,0){\line(1,0){10}}
\put(20,5){\line(1,0){10}}
\put(20,0){\line(0,1){5}}
\put(30,0){\line(0,1){5}}

\put(25,12.5){\makebox(0,0)[cc]{$G$}}
\put(25,2.5){\makebox(0,0)[cc]{$\Delta$}}

\put(3,14){\makebox(0,0)[cc]{$g$}}
\put(15,14){\makebox(0,0)[cc]{$w$}}
\put(47,4){\makebox(0,0)[cc]{$f$}}
\put(35,4){\makebox(0,0)[cc]{$v$}}
%\put(15,4){\makebox(0,0)[cc]{$h$}}
%\put(35,14){\makebox(0,0)[cc]{$y$}}
\end{picture}
	\caption{The Lur'e system}
	\label{fig:lure}
\end{figure}

A classical topic in control theory is the Lur'e problem\cite{Lurie:1944}, which concerns the stability of a feedback interconnection between a linear time-invariant (LTI) system $G$ and any nonlinearity or uncertainty $\Delta$ within some  classes (see Fig.~\ref{fig:lure}\footnote{The Lur'e system or Lur'e problem is originally defined for unforced systems, but it is relaxed to forced systems \cite{Altshuller:2013}.}). The stability of Lur'e systems is mainly studied with two different techniques:  Lyapnov stability and input-output stability.  Lyapnov stability is based on internal state variables of the unforced system, while  input-output stability studies the input-output mapping of the forced system. These two methods are closely related, and sometimes equivalent \cite{Vidyasagar:2002, Hassan:2002}. 

%The nput-output stability approach is often analysed in the frequency domain with the multiplier technique, where the class of Zames-Falb multipliers $\mathcal{M}$, defined in \cite{Zames:1966,Zames:1968} (also refer to \cite{Joaquin:2016} for a tutorial), is the widest class of LTI multipliers, and some other multipliers are phase equivalent to corresponding Zames-Falb multipliers \cite{Joaquin:2013}. A classical problem in the Zames-Falb theorem is to study   $\ell_2$-stability of the feedback interconnection in Fig. \ref{fig:lure} with a memoryless nonlinearity with the slope restricted in the range $[0,K]$. Assume the Lur'e system is in discrete time. The well-known condition for this problem is to search for a multiplier  $M\in \mathcal{M}$, such that 

The input-output approach splits the stability problem into two steps. Firstly, we need to find a class of LTI systems referred to as multipliers, preserving some properties of the set of nonlinearities $\Delta$. Secondly, we search for a suitable multiplier within the developed class of multipliers for the LTI system $G$. As a result, the stability of the nonlinear system is translated into an LTI design problem. The nonlinearity $\Delta$ is a slope-restricted nonlinearity, where the class of Zames-Falb multipliers $\mathcal{M}$, defined in both continuous-time~\cite{Zames:1968} and discrete-time~\cite{Jan:1971} (see \cite{Joaquin:2016} for a tutorial), is the widest class of LTI multipliers preserving the positivity of the nonlinearity. Note that some other multipliers are phase equivalent to corresponding Zames-Falb multipliers \cite{Joaquin:2013}. Then the absolute stability problem is reduced to a search of $M\in\mathcal{M}$ such that 
\begin{equation}\label{eq:re_condition}
Re\{M(z)\left(1+KG(z)\right)\}<0 \quad    \forall |z|=1,
\end{equation}
where $K$ is the maximum slope of the nonlinearity. This condition can be expressed equivalently in the IQC framework \cite{Alexandre:1997},  where the frequency domain condition can be converted to computable linear matrix inequalities (LMIs) by the Kalman-Yakubovich-Popov (KYP) lemma \cite{Rantzer:1996}. % Before that, a matrix related to the dynamic multiplier should be factorised. 
In discrete-time, the search over FIR Zames-Falb multipliers proposed in~\cite{Shuai:2014,Joaquin:2018} provides the less conservative absolute stability results in the literature.

Recently, the analysis of the convergence rates of the Lur'e system has attracted much attention. First-order optimisation algorithms, such as gradient decent method and Nesterov method, are written as Lur'e systems \cite{Lessard:2016}.  Specially, a strongly convex function converges to the optimal point exponentially with a rate $\rho$ ($0<\rho<1$) with first-order optimisation algorithms, which can be considered as the equilibrium point of the corresponding Lur'e system. Less-conservative results are obtained in optimisation \cite{Scoy:2018, Cyrus:2018, Mahyar:2017} and control \cite{Zachary:2017}. 

%In the $\rho$-IQC  analysis, the estimation of a tight bound of the system state exponential convergence rates in a stable Lur'e system is a main task. In the existing literature \cite{Boczar:2015,Boczar:2017,Bin:2016},   exponential stability of the system in Fig. \ref{fig:lure} is proved to be equivalent to  $\ell_2$-stability of the system in Fig. \ref{fig:lure_scaled}. Based on this relation, there are two approaches to study stability. 

\begin{figure}[b!]
	\centering
	\ifx\JPicScale\undefined\def\JPicScale{1}\fi
\unitlength \JPicScale mm
\begin{picture}(80,25)(0,0)

\linethickness{0.3mm}
\put(10,7.5){\line(1,0){10}}
\put(10,7.5){\line(0,1){7.5}}
\put(10,15){\vector(0,1){0.12}}
\put(10,17.5){\circle{5}}

\put(0,17.5){\line(1,0){7.5}}
\put(7.5,17.5){\vector(1,0){0.12}}
\put(12.5,17.5){\line(1,0){7.5}}
\put(20,17.5){\vector(1,0){0.12}}

\put(20,15){\line(1,0){10}}
\put(20,20){\line(1,0){10}}
\put(20,15){\line(0,1){5}}
\put(30,15){\line(0,1){5}}

\put(20,5){\line(1,0){10}}
\put(20,10){\line(1,0){10}}
\put(20,5){\line(0,1){5}}
\put(30,5){\line(0,1){5}}

\put(30,17.5){\line(1,0){5}}
\put(30,7.5){\line(1,0){5}}

\put(35,15){\line(1,0){10}}
\put(35,20){\line(1,0){10}}
\put(35,15){\line(0,1){5}}
\put(45,15){\line(0,1){5}}

\put(35,5){\line(1,0){10}}
\put(35,10){\line(1,0){10}}
\put(35,5){\line(0,1){5}}
\put(45,5){\line(0,1){5}}

\put(25,17.5){\makebox(0,0)[cc]{$\rho_+$}}
\put(25,7.5){\makebox(0,0)[cc]{$\rho_-$}}

\put(30,17.5){\line(1,0){5}}
\put(35,17.5){\vector(1,0){0.12}}
\put(30,7.5){\line(1,0){5}}
\put(30,7.5){\vector(-1,0){0.12}}

\put(35,15){\line(1,0){10}}
\put(35,20){\line(1,0){10}}
\put(35,15){\line(0,1){5}}
\put(45,15){\line(0,1){5}}

\put(35,5){\line(1,0){10}}
\put(35,5){\line(1,0){10}}
\put(35,5){\line(0,1){5}}
\put(45,5){\line(0,1){5}}

\put(40,17.5){\makebox(0,0)[cc]{$G$}}
\put(40,7.5){\makebox(0,0)[cc]{$\Delta$}}

\put(45,17.5){\line(1,0){5}}
\put(50,17.5){\vector(1,0){0.12}}
\put(45,7.5){\line(1,0){5}}
\put(45,7.5){\vector(-1,0){0.12}}

\put(50,15){\line(1,0){10}}
\put(50,20){\line(1,0){10}}
\put(50,15){\line(0,1){5}}
\put(60,15){\line(0,1){5}}

\put(50,5){\line(1,0){10}}
\put(50,10){\line(1,0){10}}
\put(50,5){\line(0,1){5}}
\put(60,5){\line(0,1){5}}

\put(55,7.5){\makebox(0,0)[cc]{$\rho_+$}}
\put(55,17.5){\makebox(0,0)[cc]{$\rho_-$}}

\put(60,17.5){\line(1,0){10}}
\put(70,17.5){\line(0,-1){7.5}}
\put(70,10){\vector(0,-1){0.12}}
\put(70,7.5){\circle{5}}
\put(60,7.5){\line(1,0){7.5}}
\put(60,7.5){\vector(-1,0){0.12}}

\put(72.5,7.5){\line(1,0){7.5}}
\put(72.5,7.5){\vector(-1,0){0.12}}

\put(3,19){\makebox(0,0)[cc]{$g$}}
\put(77,9){\makebox(0,0)[cc]{$f$}}

\put(15,19){\makebox(0,0)[cc]{$w$}}
\put(65,9){\makebox(0,0)[cc]{$v$}}

\put(32.5,19){\makebox(0,0)[cc]{$u_1$}}
\put(47.5,19){\makebox(0,0)[cc]{$y_1$}}

\put(32.5,9){\makebox(0,0)[cc]{$u_2$}}
\put(47.5,9){\makebox(0,0)[cc]{$y_2$}}

\put(40,22){\makebox(0,0)[cc]{$\overbrace{\quad\quad\quad\quad\quad\quad\quad\quad\quad\quad\quad\quad}$}}
\put(40,25){\makebox(0,0)[cc]{$G_{\rho}$}}

\put(40,3){\makebox(0,0)[cc]{$\underbrace{\quad\quad\quad\quad\quad\quad\quad\quad\quad\quad\quad\quad} $}}
\put(40,0){\makebox(0,0)[cc]{$\Delta_{\rho}$}}

\end{picture}
	\caption{The scaled system in \cite{Bin:2016}}
	\label{fig:lure_scaled}
\end{figure}
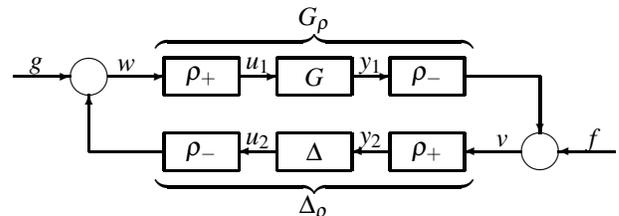
The convergence analysis of Lur'e systems has been presented in two different but equivalent frameworks:
\begin{itemize}
	\item On the one hand, in \cite{Bin:2016}, the (time domain) IQCs are constructed for the scaled uncertainty $\Delta_{\rho}$, and the corresponding exponential stability condition is to search for a multiplier that belongs to a suitable subset of $\mathcal{M}$, such that
	\begin{equation}\label{eq:re_p_condition1}
	Re\{M(z)(1+KG_{\rho}(z))\}<0 \quad \forall |z|=1,
	\end{equation} 
	where the multiplier is in the same form of original FIR Zames-Falb multipliers, while the $\ell_1$-norm condition is penalised with $\rho$.
	\item On the other hand, in \cite{Boczar:2015, Boczar:2017, Freeman:2018}, the (frequency domain)  $\rho$-IQCs are constructed for the original uncertainty $\Delta$, and the exponential stability condition is to search for a multiplier that also belongs to a subset of $\mathcal{M}$, such that
	\begin{equation}\label{eq:re_p_condition}
	Re\{M(\rho z)(1+KG(\rho z))\}<0 \quad \forall |z|=1,
	\end{equation}
	where the multiplier is constructed from the original FIR Zames-Falb multiplier $M(z)$ by replacing $z$ by $\rho z$,  and the $\ell_1$-norm condition  is  also penalised with  $\rho$.
\end{itemize}

In both approaches, sound analysis on the causal FIR Zames-Falb multipliers are provided in the literature above.  On contrast, the noncausal FIR Zames-Falb multiplier in the form $M(\rho z)$ is studied in \cite{Freeman:2018}, where its modified $\ell_1$-norm condition is proved, but  details to obtain the stability LMI are not given.

In this technical note, we are concerned with the technique to apply noncausal FIR Zames-Falb multipliers in \cite{Shuai:2014,Joaquin:2018} to estimate exponential convergence rates, and especially focus on the factorisations. The main contribution of this technical note is the development of suitable factorizations for both approaches when noncausal Zames-Falb multipliers are used. %In the literature above, a standard factorisation proposed in \cite{Heath:2005} is used for causal multipliers, where the direct replacement from $z$ to $\rho z$ is valid in the $\rho$-IQC framework. However, this factorisation is invalid for anticausal or noncausal multipliers.  Besides, we will show that not every factorisation for noncausal Zames-Falb IQC can be applied to the $\rho$-IQC case with the multiplier $M(\rho z)$. 
In Section \ref{sec:factorisations_iqc}, the time domain Zames-Falb IQC with causal multipliers in \cite{Bin:2016} for continuous time system is extended to frequency domain with noncausal multipliers for  discrete time system. 
Meanwhile, in Section \ref{sec:factorisations_piqc}, we provide a factorisation by lifting \cite{Hosoe:2013} as an unified structure for causal, anticausal and noncausal FIR Zames-Falb multipliers in both the IQC and $\rho$-IQC frameworks. Then, the validity of different factorisations are discussed, which completes the results in \cite{Boczar:2015,Boczar:2017,Freeman:2018}.  Furthermore, we show that the multipliers $M(z)$ in (\ref{eq:re_p_condition1}) and $M(\rho z)$ in (\ref{eq:re_p_condition})  are equivalent by variable conversion, and lead to similar results in numerical examples in Section \ref{sec:results}.

\section{Notations and preliminary results}\label{sec:notations}
Some of the notations and definitions are summarised from \cite{Lessard:2016, Boczar:2015, Boczar:2017, Bin:2016}, which are repeated here for completeness.

\subsection{Notations and Lur'e systems}

Let $\mathbb{Z}$ and $\mathbb{Z}^+$ be the set of integer numbers and positive integer numbers including zero, respectively. The notations $\mathbb{R}$ and $\mathbb{R}^+$ are defined in the same way for real numbers. And let $\mathbb{C}$ be the set of complex numbers.  Let $\ell$ be the space of all real-valued sequences $h: \mathbb{Z}^+\mapsto\mathbb{R}$. Let $\ell_2$ be the space of real-valued square-summable sequences $h:\mathbb{Z}^+\mapsto\mathbb{R}$. For any $\rho\in[0,1]$, we will say that $h\in\ell_2^\rho$ if the sequence $\{\rho^{-k}h_k\}_{k=0}^\infty$ belongs to $\ell_2$. Finally, for absolute-summable sequences $h: \mathbb{Z}^+\mapsto\mathbb{R}$, we define $\|h\|_1=\sum_{k=-\infty}^{\infty}|h_k|$.

Let $\mathbf{RH}_{\infty}$  and $\mathbf{RL}_{\infty}$ be the space consisting of proper real rational transfer functions $G$:  $G\in \mathbf{RH}_{\infty}$ has all poles inside the open unit disk in the complex plane;   $G\in \mathbf{RL}_{\infty}$ has no pole on the unit disk. With the minimal state-space realisation, the transfer function is $G(z) = C(zI -A)^{-1}B+D$, or $G \sim $$\begin{bmatrix}  A\;  B;\; C\; D \end{bmatrix}$ in short. 
The expression $G^{*}(z)$ denotes the complex conjugate transpose of $G(z)$ at $|z|=1$, i.e. $G^{*}(z)={G}^T\left(\frac{1}{z}\right)$, 
where the superscript $T$ indicates the transpose. Moreover, if a parameter $\rho$ is involved in the variable, the complex conjugate transpose can be expressed as $G^*(\rho, z)=G^T\left(\rho, \frac{1}{z}\right)$.

A nonlinear operator $\Delta: \ell(\mathbb{Z}^+) \mapsto \ell(\mathbb{Z}^+)$ is said to be memoryless if there exists a map $N:  \mathbb{R} \to \mathbb{R} $ such that $(\Delta\upsilon)_k=N(\upsilon_k)$, $\forall k \in \mathbb{Z}$. Assume that $\Delta(0)=0$. The memoryless uncertainty $\Delta$ is said to be (sector) bounded, denoted by $\Delta \in [\underline{k},\overline{k}]$ ($0\le \underline{k}<\overline{k}<\infty$), if $\underline{k} x \le N(x) \le \overline{k} x, \forall x \in \mathbb{R}$. The uncertainty $\Delta$ is said to be slope-restricted, denoted by $\Delta\in S[\underline{k},\overline{k}]$, if $\underline{k}(x_1-x_2) \le N(x_1)-N(x_2) \le \overline{k}(x_1-x_2), \forall x_1, x_2 \in \mathbb{R}$ and $x_1 \ne x_2$. The slope-restricted uncertainty is also sector bounded, but the reverse is not.  Finally, the  uncertainty $\Delta$ is said to be odd if $\Delta(-x)=-\Delta(x)$, $\forall x\in \mathbb{R}$.

Consider the Lur'e system in Figure \ref{fig:lure}. It is expressed as
\begin{equation*}
v=f+Gw, \quad w=g+\Delta v.
\end{equation*}

The feedback interconnection is well-posed if the inverse map $(v,w) \mapsto (g,f)$ is causal in $\ell$.

\begin{definition}
The feedback  interconnection  in Fig. \ref{fig:lure} is  $\ell_2$-stable if it is well-posed, and the signals $(v,w)\in\ell_2$ for any $(g,f)\in\ell_2$.
\end{definition} 

%$\ell_2$-stable with finite gain if it is well-posed and there exists a constant $\gamma$ such that  
%\begin{equation}\label{eq:bibo_stable}
%\|w\|+\|v\|<\gamma(\|g\|+\|f\|).
%\end{equation}
%\end{definition} 
%\begin{remark}
%Condition (\ref{eq:bibo_stable}) implies  $\ell_2$-stability: the signals $(v,w)$ are in $\ell_2$ for any $(g,f)$ in $\ell_2$.  
%\end{remark}

\begin{definition} The feedback  interconnection in Fig. \ref{fig:lure} is  globally exponentially stable with convergence rate $\rho$ if there exists some $\rho\in (0,1)$ and $c>0$ such that when $g=0$ and $f=0$,
\begin{equation}\label{eq:exponential_stable}
\|x_k\|\le c\rho^k \|x_0\| \quad \forall k\ge 0,\;\forall x_0\in \mathbb{R}^n .
\end{equation}
\end{definition} 

Henceforth, the infimum convergence rate of the feedback  interconnection in Fig. \ref{fig:lure} is referred to as $\rho_{\{G,\Delta\}}$. 

\begin{remark}\label{rm:exponential_convergence}
Condition (\ref{eq:exponential_stable}) is  equivalent  to the fact that  the state $x_k$ of $G$ converges to zero exponentially with the rate $\rho$, i.e. $x_k \rho^{-k}\to 0$ as $k\to \infty$.
\end{remark}

As mentioned in the introduction,  $\ell_2$-stability is an input-output relation, while  exponential stability is an internal relation. Therefore, it is not trivial to restate  exponential stability in an input-output manner.  

\begin{definition}\label{df:l2p_stability}
The feedback  interconnection  in Fig. \ref{fig:lure} is  $\ell_2^{\rho}$-stable if it is well-posed, and the signals $(v,w)\in\ell_2^{\rho}$ for any $(g,f)\in\ell_2^{\rho}$.
\end{definition}

\begin{theorem}\label{th:stability_relations}
For the Lur'e system in Fig. \ref{fig:lure}, assume $G$ is controllable and observable, and $\Delta$ is  memoryless and slope-restricted. The unforced system is globally exponentially stable with rate $\rho$ if and only if  the forced system is  $\ell_2^{\rho}$-stable.
\end{theorem}

\begin{IEEEproof}
The sufficiency can be proved in a similar way with Proposition 5 in \cite{Boczar:2015,Boczar:2017}, and the necessity  is  proved in outline in Appendix \ref{appendix_proof}. 	
\end{IEEEproof}

%Notice that the exponential stability is equivalent to the BIBO stability when $\rho=1$.

\subsection{Kalman conjecture for convergence analysis}\label{sc:kalman_conjecture}
The Kalman conjecture is a necessary and sufficient condition for stability when it is true, which is stated below.

\begin{definition}[Nyquist value, $K_N$] \label{df:Nyquist}
	The Nyquist value of a stable transfer function $G(z)$ is
	\begin{equation*}\label{Nyquist}
	K_N=\sup_{K}\{ K>0: (1-\tau KG(z))^{-1} \ \textrm{is stable} \ \forall \tau\in[0,1] \}.
	\end{equation*}    
\end{definition}

\begin{conjecture}[Kalman conjecture~\cite{Kalman:1957}] Let $\Delta$ be memoryless,  and $\Delta \in S[0, K]$. The feedback interconnection between $G$ and $\Delta$  is asymptotically stable if and only if  $K < K_N$.
\end{conjecture}

We can translate the above definition and conjecture into the convergence analysis. Let us define the absolute convergence rate of the class of systems defined by the Lur'e system as follows
\begin{equation}\label{eq:rho_star}
\rho^*_{\{G,K\}}=\sup_{\Delta\in S[0, K]}\{\rho_{\{G,\Delta\}}\}
\end{equation}
where $K<K_N$.
A lower bound of the $\rho^*_{\{G,K\}}$ is given by
\begin{equation}\label{eq:rho_linear}
\underline{\rho^*_{\{G,K\}}}=\max_{\tau\in[0,1]}\left\{\left|\text{eig}\left(\frac{G}{1-\tau KG}\right)\right|\right\}.
\end{equation}
In some instances, this lower bound is referred to as the theoretical value.  The Kalman conjecture can be restated using the convergence rates defined in (\ref{eq:rho_star}) and (\ref{eq:rho_linear}) as follows.
\begin{conjecture}[Kalman conjecture for convergence analysis]\label{cj:kc_convergence}
	 For any stable $G$, let $\Delta \in S[0, K]$ with $K<K_N$, then
\begin{equation}\label{eq:1}
\rho^*_{\{G,K\}}=\underline{\rho^*_{\{G,K\}}}.
\end{equation}	
\end{conjecture}
%{\color{red} (What is the accurate definition of $\rho^*$?) Similar to the analysis in the Kalman conjecture, the theoretical  value $\rho^*$  of  the lower bound of the worst-case convergence rates is given by the linear feedback case,}

%where $\tau \in [0,1]$ and $K < K_N$. The value $\rho^*$ may not be reached when the Kalman conjecture is wrong. In this case, sufficient conditions for the stability should be used to estimate a tight upper bound of $\rho$. 

\subsection{Estimation of upper bound of $\rho^*_{\{G,K\}}$}\label{sc:estimation}
In the last years, Zames-Falb multipliers have been used to estimate an upper bound of $\rho^*_{\{G,K\}}$, denoted by $\overline{\rho^*_{\{G,K\}}}$. As mentioned in the introduction, there are two approaches based on the relation below.

\begin{theorem}[\cite{Boczar:2015,Bin:2016,Boczar:2017}]\label{th:stability_original_sclaed}
		The system in Fig. \ref{fig:lure} is well-posed if and only if the scaled system in Fig. \ref{fig:lure_scaled} is well-posed. Furthermore, the system in Fig. \ref{fig:lure} is $\ell_2^{\rho}$-stable if and only if the scaled system in Fig. \ref{fig:lure_scaled} is $\ell_2$-stable. 
\end{theorem} 

%\begin{IEEEproof}
%In Fig. \ref{fig:lure_scaled}, the operators $\rho_+$ and $\rho_-$ are defined in time domain as $\rho^k$ and $\rho^{-k}$, respectively. Then, $\rho_-\circ(G(z)\circ \rho_+) \equiv G(\rho z)$, $\rho_-\circ(\Delta\circ \rho_+) \equiv \Delta_{\rho}$.
%Moreover,  $v(z) \equiv y_2(\rho z)$, $(\Delta_{\rho}(v))(z) \equiv u_2(\rho z)$, and similar relations hold for $G_{\rho}$. As a result, the signals $(w, G_{\rho}(w))$ and $(v, \Delta_{\rho}(v)) \in \ell_2$ if and only if $(u_1, y_1)$ and $(y_2, u_2)\in \ell_2^{\rho}$, respectively. See the formal proof in \cite{Boczar:2017}, which is not repeated here.
%\end{IEEEproof}

The two approaches to estimate  $\overline{\rho^*_{\{G,K\}}}$ are reviewed in the following parts.

\subsubsection{Analysis in IQC framework}

In this approach, $\ell_2^{\rho}$ stability of the system in Fig. \ref{fig:lure} is studied by  $\ell_2$ stability of the scaled system in Fig. \ref{fig:lure_scaled}, where an IQC is constructed for the scaled uncertainty $\Delta_{\rho}$ at first.

\begin{definition}[IQC \cite{Alexandre:1997}]\label{df:IQC}	
	Let $\Pi(z)$  be a Hermitian (self-adjoint) bounded measurable operator.  Then, for a bounded and causal operator $\Delta_{\rho}: {\ell} \mapsto {\ell}$, it is said to satisfy the IQC defined by $\Pi$, if  for all $v\in {\ell}_2$
	\begin{gather}\label{eq:iqc}
	\int_{|z|=1}
	\begin{bmatrix}
	 \hat{v}(z) \\ \widehat{\Delta_{\rho}v}(z)
	\end{bmatrix}^*
	\Pi(z)
	\begin{bmatrix}
	 \hat{v}(z) \\ \widehat{\Delta_{\rho}v}(z)
	\end{bmatrix}
	dz \ge 0,
	\end{gather}
where $\hat{v}$ and  $\widehat{\Delta_{\rho}v}$ denote the z-transform of $v$ and $\Delta_{\rho}v$ respectively. 
\end{definition}

\begin{theorem}[\cite{Alexandre:1997}] \label{th:IQC_discrete}
	For the system in Fig. \ref{fig:lure_scaled}, let $G_{\rho}(z) \equiv G(\rho z) \in \mathbf{RH}_{\infty}$,  and $\Delta_{\rho}$ be a causal bounded operator. Assume that $\forall \tau \in [0,1]$,
	\begin{enumerate}
		\item  the feedback interconnection between $G_{\rho}$ and $\tau \Delta_{\rho}$ is well-posed;
		\item the operator $\tau \Delta_{\rho}$ satisfies the IQC defined by $\Pi$;
		\item there  exists $\epsilon>0$, such that
		\begin{gather}\label{eq:iqc_stability}
		\begin{bmatrix}
		G_{\rho}(z)  \\
		I \\
		\end{bmatrix}^*
		\Pi(z)
		\begin{bmatrix}
		G_{\rho}(z)  \\
		I \\
		\end{bmatrix}
		\le -\epsilon I, \quad \forall |z|=1.
		\end{gather}
	\end{enumerate}
 Then, the system in Fig.\ref{fig:lure_scaled} is $\ell_2$-stable, thus the system in Fig.\ref{fig:lure} is $\ell_2^{\rho}$-stable. 
\end{theorem}

In order to make the frequency domain inequality (FDI)  (\ref{eq:iqc_stability}) computable, the Kalman-Yakubovich-Popov (KYP) lemma should be applied.

\begin{lemma}[KYP lemma \cite{Rantzer:1996}]\label{Le:kyp}
	Given $A \in \mathbb{R}^{n \times n}$, $B \in \mathbb{R}^{n \times m}$, $K_p=K_p^{T} \in  \mathbb{R}^{(n+m) \times (n+m)}$, with $ det(z I -A)\ne 0$ for all $|z|=1$, where the pair $(A, B)$ are controllable, the following statements are equivalent:
	\begin{enumerate}
		\item For all $|z|=1$,
		\begin{gather*}
		\begin{bmatrix}
		(zI -A)^{-1}B  \\
		I \\
		\end{bmatrix}^*
		K_p
		\begin{bmatrix}
		(zI -A)^{-1}B  \\
		I \\
		\end{bmatrix}
		\leq 0.
		\end{gather*}
		\item There is a symmetric matrix $P \in \mathbb{R}^{n \times n}$ and
		\begin{gather*}
		\begin{bmatrix}
		A^{T}PA-P & A^TPB   \\
		B^{T}PA & B^TPB \\
		\end{bmatrix}+K_p
		\leq 0.
		\end{gather*}
	\end{enumerate}
\end{lemma}

Generally, the IQC multiplier $\Pi$ is dynamic, and can be factorised as below.
\begin{definition}[\cite{Scherer:2011}]
	Any $\Pi(z)\in\mathbf{RL}_\infty$ has nonunique factorisations $(\Psi,K_p)$ in the form 
	\begin{equation}\label{eq:factorisation}
	\Pi(z)=\Psi^*(z)K_p\Psi(z),
	\end{equation}
	where $K_p=K_p^T$ is constant, and $\Psi$ is a stable LTI system with the state-space representation
	\begin{gather}\label{eq:psi_state_space}
	\Psi(z) \sim  \begin{bmatrix}  A_\Psi & B_{\Psi_1} & B_{\Psi_2}\\ C_\Psi & D_{\Psi_1} & D_{\Psi_2} \end{bmatrix}.
	\end{gather}
\end{definition}

Notice that when $G(z) \sim $ $\begin{bmatrix}  A\;  B;\; C\; D \end{bmatrix}$, $G_{\rho}(z)\equiv $ $G(\rho z) $ $\sim \begin{bmatrix}  \rho^{-1}A\;\;  \rho^{-1}B;\;\;  C\;\;  D \end{bmatrix}$. Next, substitute (\ref{eq:factorisation}) into (\ref{eq:iqc_stability}), and apply the KYP lemma, the well-known stability LMI is applied to the scaled system.

\begin{corollary} 
	The FDI (\ref{eq:iqc_stability}) is equivalent to the existence of $P=P^T$ such that
	\begin{equation}\label{eq:iqc_lmi}
	\begin{bmatrix} \hat{A}^TP\hat{A}-P  &\hat{A}^TP\hat{B} \\  \hat{B}^TP\hat{A} & \hat{B}^TP\hat{B}\end{bmatrix}+\begin{bmatrix} \hat{C}^T \\\hat{D}^T\end{bmatrix} K_p \begin{bmatrix} \hat{C} &\hat{D}\end{bmatrix}<0,
	\end{equation}
	where $\Psi\begin{bmatrix} G_{\rho}  \\I\end{bmatrix} \sim \begin{bmatrix} \hat{A} & \hat{B}\\ \hat{C} & \hat{D} \end{bmatrix}$, and
	$\hat{A}=\begin{bmatrix}\rho^{-1}A & 0 \\ B_{\Psi_1}C & A_{\Psi} \end{bmatrix}$,
	$\hat{B}=\begin{bmatrix}\rho^{-1}B \\ B_{\Psi_2}+B_{\Psi_1}D \end{bmatrix}$,
	$\hat{C}=\begin{bmatrix}D_{\Psi_1}C & C_{\Psi}\end{bmatrix}$, $ \hat{D}=D_{\Psi_2}+D_{\Psi_1}D$.
\end{corollary}

\subsubsection{Analysis in $\rho$-IQC framework}

In this approach,  $\ell_2^{\rho}$ stability of the system in Fig. \ref{fig:lure} is studied by scaling signals,  where a $\rho$-IQC is defined for the original uncertainty $\Delta$ at first.

\begin{definition}[$\rho$-IQC \cite{Boczar:2015,Boczar:2017}] \label{df:rho_IQC}
	Let $\Pi(\rho, z)$ be a Hermitian (self-adjoint, i.e. $\Pi(\rho, z)=\Pi^*(\rho, z)=\Pi^T(\rho, z^{-1})$) bounded measurable operator.  Then, for a bounded and causal operator $\Delta: {\ell} \mapsto {\ell}$, it is said to satisfy the $\rho$-IQC defined by $\Pi$, if for all $y_2\in {\ell}_2^{\rho}$
			\begin{gather}\label{eq:rho_iqc}
			\int_{|z|=1}
			\begin{bmatrix}
			\widehat{y_2}(\rho, z)\\ \widehat{\Delta y_2}(\rho, z)
			\end{bmatrix}^*
			\Pi(\rho, z)
			\begin{bmatrix}
			\widehat{y_2}(\rho, z)\\ \widehat{\Delta y_2}(\rho, z)
			\end{bmatrix}
			dz \ge 0,
			\end{gather}
	where $\widehat{y_2}(\rho, z) \equiv \widehat{y_2}(\rho z)$ and $\widehat{\Delta y_2}(\rho, z)\equiv \widehat{\Delta y_2}(\rho z)$. 
\end{definition}

\begin{remark}
According to the signal relations in the proof of Theorem \ref{th:stability_original_sclaed}, (\ref{eq:iqc}) holds if and only if  (\ref{eq:rho_iqc}) holds, which implies $\Pi(z)$ and $\Pi(\rho,z)$ are  equivalent. It will be shown by the Zames-Falb IQC and $\rho$-IQC in Section \ref{sec:factorisations_iqc} and \ref{sec:factorisations_piqc}.
\end{remark}

\begin{theorem} [\cite{Boczar:2015,Boczar:2017}] \label{th:rho_IQC}
	Fix $\rho\in(0,1)$. For the Lur'e system in Fig. \ref{fig:lure}, let $G(\rho,z) \equiv G(\rho z) \in \mathbf{RH}_{\infty}$, and  $\rho_-\circ(\Delta\circ \rho_+) \equiv \Delta_{\rho}$ be a causal bounded operator. Assume that $\forall \tau \in [0,1]$,
	\begin{enumerate}
		\item  the feedback interconnection between $G$ and $\tau \Delta$ is well-posed;
		\item the operator $\tau \Delta$ satisfies the $\rho$-IQC defined by $\Pi$;
		\item there  exists $\epsilon>0$, such that
		\begin{gather}\label{eq:rho_iqc_stability}
		\begin{bmatrix}
		G(\rho, z)  \\
		I \\
		\end{bmatrix}^*
		\Pi(\rho, z)
		\begin{bmatrix}
		G(\rho, z)  \\
		I \\
		\end{bmatrix}
		\le -\epsilon I, \forall |z|=1.
		\end{gather}
	\end{enumerate}
 Then, the system in Fig.\ref{fig:lure} is $\ell_2^{\rho}$-stable.
\end{theorem}

Similarly, (\ref{eq:rho_iqc_stability}) can be converted to the stability LMI after the factorisation of $\Pi(\rho,z)$ as defined below.

\begin{definition}
	Any $\Pi(\rho,z) \in\mathbf{RL}_\infty$ has nonunique factorisations $(\Psi,K_p)$ in the form 
	\begin{equation}\label{eq:factorisation_rho}
	\Pi(\rho,z)=\Psi^*(\rho,z)K_p\Psi(\rho,z),
	\end{equation}
	where $K_p=K_p^T$ is constant, and $\Psi(\rho,z)$ is a stable LTI system with the variable $z$ and the parameter $\rho$.
\end{definition}

Particularly, when the multiplier used in the $\rho$-IQC is causal, $\Psi(\rho,z)=\Psi(\rho z)$  is valid. Additionally, when (\ref{eq:psi_state_space}) holds, 
\begin{equation}\label{eq:psi_pz_causal}
\Psi(\rho, z) \sim  \begin{bmatrix}  \rho^{-1}A_\Psi & \rho^{-1}B_{\Psi_1} & \rho^{-1}B_{\Psi_2}\\ C_\Psi & D_{\Psi_1} & D_{\Psi_2} \end{bmatrix}. 
\end{equation}

Therefore, $\Psi(\rho, z)\begin{bmatrix} G(\rho, z)  \\I\end{bmatrix}$ $\sim$ $\begin{bmatrix} \rho^{-1}\hat{A} & \rho^{-1}\hat{B}\\ \hat{C} & \hat{D} \end{bmatrix}$, and the LMI   can be further simplified to the form in \cite{Boczar:2015, Lessard:2016, Boczar:2017}: $\exists P=P^T$, such that
\begin{equation}\label{eq:causal_LMI}
\begin{bmatrix} \hat{A}^TP\hat{A}-\rho^2 P  &\hat{A}^TP\hat{B} \\  \hat{B}^TP\hat{A} & \hat{B}^TP\hat{B}\end{bmatrix}+\begin{bmatrix} \hat{C}^T \\\hat{D}^T\end{bmatrix} K_{p} \begin{bmatrix} \hat{C} &\hat{D}\end{bmatrix} <0.
\end{equation}

\begin{figure}[ht]
\centering
\includegraphics[width=\linewidth]{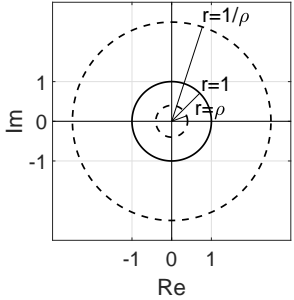}
\caption{Analysis on the circles with radii $\rho$ and $1/\rho$ }
\label{fig:circles}
\end{figure}

	\begin{remark}
		In the conventional IQC framework, the variable is $z$ only, and the analysis is conducted on the circle $|z|=1$ in the complex plane. However, in the $\rho$-IQC framework, the variables are $\rho$ and $z$, such as $\rho z$ and $\frac{z}{\rho}$. Then, when the complex conjugate  is used, the analysis is conducted on both circles $|z|=\rho$ and $|z|=\frac{1}{\rho}$, as shown in Fig. \ref{fig:circles}.
	\end{remark}

\subsection{Zames-Falb IQC with FIR multipliers} \label{sc:ZF-IQC}
In this part, the structure of the Zames-Falb IQC for the class of slope-restricted uncertainties is introduced. Then, three factorisations with FIR multipliers are provided.

%\begin{theorem}[Zames-Falb IQC]
%	Assume the uncertainty $\Delta$ is static and $\Delta\in S[0,\infty)$. It satisfies the Zames-Falb IQC defined by $\Pi_0$ as
%	\begin{equation}\label{eq:zf_iqc1}
%	\Pi_0(z)=\begin{bmatrix}0 & M^*(z) \\ M(z)  & 0 \end{bmatrix}, 
%	\end{equation}
%	where $M(z)$ is a Zames-Falb multiplier. 
%\end{theorem}  

\begin{theorem}[Zames-Falb IQC \cite{Heath:2005}]
Assume the uncertainty $\Delta$ is static and $\Delta\in S[0, K]$. It satisfies the Zames-Falb IQC defined by $\Pi$ as
\begin{equation}\label{eq:zf_iqc}
%\begin{split}
\Pi(z)=
% \begin{bmatrix}
%\beta & -1\\ -\alpha &1
%\end{bmatrix}^T
%\Pi_0(z)
%\begin{bmatrix}
%\beta & -1\\ -\alpha &1
%\end{bmatrix},\\
\begin{bmatrix} 0 & K M^*(z) \\ K M(z) & -(M(z)+M^*(z)) \end{bmatrix}.
%\end{split}
\end{equation}
\end{theorem}

Here, we focus on the noncausal FIR Zames-Falb multiplier proposed in \cite{Shuai:2014},
\begin{equation}\label{eq:zf_multiplier}
M(z)=-h_{-n_b}z^{-n_b}-\cdots-h_{-1}z^{-1}+h_0-h_{1}z^{1}-\cdots-h_{n_f}z^{n_f},
\end{equation} 
where the causal part is with the backward-shift operator $z^{-i_b}$ ($i_b=1,2, \cdots, n_b$), and the anticausal part is with the forward-shift operator $z^{i_f}$ ($i_f=1,2, \cdots, n_f$). In addition, $h_{i_f}>0$ and $h_{-i_{b}}>0$, or $\Delta$ is odd. The $\ell_1$-norm  condition of $M(z)$ is
\begin{equation}\label{eq:FIR_L1}
\sum_{i_b=1}^{n_b}|h_{-i_b}| + \sum_{i_f=1}^{n_f}|h_{i_f}|< h_0,
\end{equation}
where we can set $h_0=1$ without loss of generality.

A standard factorisation of (\ref{eq:zf_iqc}) used in previous literature, such as \cite{Boczar:2015, Lessard:2016}, is (\ref{eq:factorisation}) with 
	\begin{equation}\label{eq:zf_factorisation1}
	\Psi_1(z)=\begin{bmatrix} K M(z) & -M(z)  \\0 & 1 \end{bmatrix}, \quad K_{p,1}=\begin{bmatrix} 0 & 1  \\1 & 0 \end{bmatrix},
	\end{equation}
where $M(z)$ must be causal to keep $\Psi_1(z)$ stable, so it  needs  further factorisation for noncausal multipliers. The state-space representation of $\Psi_1(z)$ with causal Zames-Falb multipliers is given in \cite{Lessard:2016}.

Moreover,  the factorisation method called ``lifting factorisation'' is available, which can be treated as the discrete time counterpart of the factorisation for general continuous time multipliers in \cite{Veenman:2016}. One possible lifting factorisation is with
	 {\small{\begin{equation}\label{eq:zf_factorisation2}
 	\Psi_2(z)=\left[\begin{array}{c|c}
 	1 & 0\\ \pmb{Z}^{-i} & \pmb{0}\\ \hline 0 & 1\\\pmb{0} & \pmb{Z}^{-i} \end{array}\right], 
 	K_{p,2}=\arraycolsep1pt\left[\begin{array}{cc|cc} 0 & \pmb{0} & Kh_0 & K\pmb{h}^T_{i} \\ 
 	\pmb{0} & \pmb{0} & K \pmb{h}_{-i} & \pmb{0}\\ \hline 
 	Kh_0 & K\pmb{h}^T_{-i} & -2h_0 & -\pmb{h}^T_{-i}-\pmb{h}^T_{i}\\
 	K\pmb{h}_{i} & \pmb{0}  & -\pmb{h}_{-i}-\pmb{h}_{i} & \pmb{0}
 	\end{array}\right],
 	\end{equation}}}
where $\Psi_2(z)$ is called "lifting matrix", whose state-space representation is provided in Appendix \ref{appendix_ss}. Particularly, the causal and the anticausal parts in (\ref{eq:zf_multiplier}) must have the same step with this factorisation, i.e. $n_b=n_f=n_z$. Additionally, $\pmb{Z}^{-i}=[z^{-1} \; z^{-2} \;\cdots \; z^{-n_z}]^T$; $\pmb{h}_{i}=[h_1\; h_2\; \cdots \; h_{n_z}]^T$, $\pmb{h}_{-i}=[h_{-1}\; h_{-2}\; \cdots \; h_{-n_z}]^T$. Moreover, causal multipliers are obtained with $\pmb{h}_{i}=\pmb{0}$; anticausal multipliers are obtained with $\pmb{h}_{-i}=\pmb{0}$. Notice that the causal part and anticausal part share the same base $\begin{bmatrix}
1\\ \pmb{Z}^{-i} \end{bmatrix}$ in $\Psi_2(z)$, so we say they are coupled. 

Henceforth, we use the notations $\pmb{0}$ for zero matrices in some proper dimensions, and  $\pmb{I}_{(n)}$ for the $n\times n$ identity matrix.

Finally, another lifting factorisation is defined in (\ref{eq:zf_factorisation3}) on the next page.
\begin{figure*}[t!]
%	\hrule
\begin{equation}\label{eq:zf_factorisation3}
\Psi_3(z) =\left[ \begin{array}{c|c} 
1 & 0\\ \pmb{Z}^{-i_b} & \pmb{0}\\0 & 1\\\pmb{0} & \pmb{Z}^{-i_b} \\ \hline 1 & 0\\ \pmb{Z}^{-i_f} & \pmb{0}\\0 & 1\\\pmb{0} & \pmb{Z}^{-i_f}
 \end{array} \right],\quad
K_{p,3} =
\left[ \begin{array}{cccc|cccc}
0 & \pmb{0} & K{h_0/2} & \pmb{0} & 0 & \pmb{0} & 0& \pmb{0}\\
\pmb{0} & \pmb{0} & K\pmb{h}_{-i_b} & \pmb{0} & \pmb{0} & \pmb{0} & \pmb{0} & \pmb{0}\\
Kh_0/2 & K\pmb{h}_{-i_b}^{T} & -h_0 & -\pmb{h}_{-i_b}^{T} & 0 & \pmb{0} & 0 & \pmb{0}\\
\pmb{0} & \pmb{0} & -\pmb{h}_{-i_b} & \pmb{0} & \pmb{0} & \pmb{0} & \pmb{0} & \pmb{0}\\ \hline
0 & \pmb{0} & 0 & \pmb{0} & 0 & \pmb{0} &K{h_0/2} & K\pmb{h}_{i_f}^{T}\\
\pmb{0} & \pmb{0} & \pmb{0} & \pmb{0} & \pmb{0} & \pmb{0} & \pmb{0} & \pmb{0}\\
0 & \pmb{0} & 0 & \pmb{0} & K{h_0/2}  & \pmb{0} & -h_0 & -\pmb{h}_{i_f}^{T}\\
\pmb{0} & \pmb{0} & \pmb{0} & \pmb{0} & K\pmb{h}_{i_f} & \pmb{0} & -\pmb{h}_{i_f} & 0\\
\end{array} \right].
\end{equation}
	\hrule
\end{figure*}
In this factorisation, the definitions of matrices are similar with (\ref{eq:zf_factorisation2}) but with possibly different values of $n_b$ and $n_f$. Similarly,  $\pmb{h}_{i_f}=\pmb{0}$ for causal multipliers;  $\pmb{h}_{-i_b}=\pmb{0}$ for anticausal multipliers. Especially, as illustrated in the lifting matrix $\Psi_3(z)$, the bases of the causal part and anticausal part are separated, which brings the flexibility to construct asymmetric noncausal multipliers. We say the causal and anticausal parts are decoupled.

In the following sections \ref{sec:factorisations_iqc} and  \ref{sec:factorisations_piqc}, the validity of these three factorisations will be discussed in IQC and $\rho$-IQC frameworks respectively. In addition, the multipliers $M(z)$ in (\ref{eq:re_p_condition1}) and  $M(\rho z)$ in (\ref{eq:re_p_condition}) will be shown equivalent.

\section{Zames-Falb multipliers for convergence analysis in IQC framework}\label{sec:factorisations_iqc}
In the convergence analysis, the FIR Zames-Falb multipliers belong to a subset of the class of Zames-Falb multipliers ($M\in\mathcal{M}_{\rho}\subset \mathcal{M}$) as their $\ell_1$ norm conditions are penalised with the convergence rate $\rho$.

Consider the noncausal multiplier in the IQC for the scaled uncertainty $\Delta_{\rho}$,
\begin{equation}\label{eq:FIR_noncausal1}
M(z)=-\tilde{h}_{-n_b}z^{-n_b}-\cdots -\tilde{h}_{-1}z^{-1} +h_0-\tilde{h}_{1}z^{1}-\cdots-\tilde{h}_{n_f}z^{n_f},
\end{equation} 
where $\tilde{h}_{-i_b}>0$ and $\tilde{h}_{i_f}>0$,  or $\Delta_{\rho}$ is odd. Its $\ell_1$ norm condition is
\begin{equation}\label{eq:FIR_noncausal_L11}
\sum_{i_b=1}^{n_b}|\tilde{h}_{-i_b}|\rho^{-i_b}+\sum_{i_f=1}^{n_f}|\tilde{h}_{i_f}|\rho^{-i_f}<h_0,
\end{equation}
where the proof of the causal part is given in \cite{Bin:2016}, while the proof of the anticausal part will be linked with the anticausal multipliers in the next section. 

As mentioned, for causal multipliers, all the factorisations in Section \ref{sc:ZF-IQC} are valid; while for anticausal and noncausal multipliers, (\ref{eq:zf_factorisation2}) and (\ref{eq:zf_factorisation3})  are valid.  Moreover, the analysis is in the IQC-framework, so the LMI also keeps the same form in (\ref{eq:iqc_lmi}).  

In short, by this approach, everything  keeps the same as in the conventional IQC analysis except that the $\ell_1$ norm condition of FIR Zames-Falb multipliers are penalised symmetrically on causal and anticausal parts.

\section{Zames-Falb multipliers for convergence analysis in $\rho$-IQC framework}\label{sec:factorisations_piqc}

Different from the IQC analysis in the previous section, the parameter $\rho$ is involved as a variable in the $\rho$-IQC. As a result, the factorisation is restricted in different cases with causal, anticausal and noncausal multipliers.

Firstly, the lifting factorisation (\ref{eq:zf_factorisation3}) can be extended to (\ref{eq:factorisation_rho}) with $K_{p,3}$ being the same, and $\Psi_3$ being modified to 

\begin{equation}\label{eq:phi_pz}
\small{\Psi_3(\rho, z) =\left[ \begin{array}{c|c}
	1 & 0\\ \rho^{-i_b}\pmb{Z}^{-i_b} & \pmb{0}\\0 & 1\\\pmb{0} & \rho^{-i_b}\pmb{Z}^{-i_b} \\ \hline 1 & 0\\ \rho^{i_f}\pmb{Z}^{-i_f} & \pmb{0}\\0 & 1\\\pmb{0} & \rho^{i_f}\pmb{Z}^{-i_f} 
	\end{array} \right]},
\end{equation}
where $\rho^{-i_b}$ and $\rho^{i_f}$ are multiplied to $\pmb{Z}^{-i_b}$ and $\pmb{Z}^{-i_f}$, respectively. The state-space representation of $\Psi_3(\rho, z)$ is attached in Appendix \ref{appendix_ss}.

In the following, it is straightforward to show that the factorisation (\ref{eq:factorisation_rho}) with $\left(\Psi_3(\rho,z), K_{p,3}\right)$ is an unified  structure in the $\rho$-IQC framework with causal, anticausal and noncausal FIR Zames-Falb multipliers.

With the noncausal multiplier in (\ref{eq:FIR_noncausal1}), setting $h_{-i_b}=\tilde{h}_{-i_b}/\rho^{-i_b}$ and $h_{i_f}=\tilde{h}_{i_f}/\rho^{i_f}$,  the noncausal multiplier $M(\rho,z)$ in the $\rho$-IQC for the original  uncertainty $\Delta$ defined in \cite{Boczar:2015,Boczar:2017,Freeman:2018} is obtained:
\begin{equation}\label{eq:FIR_noncausal}
\begin{split}
M(\rho, z) =& -h_{-n_b}\rho^{-n_b}z^{-n_b} - \cdots - h_{-1}\rho^{-1}z^{-1}\\ & +h_0 - h_{1}\rho z - \cdots - h_{n_f}\rho^{n_f}z^{n_f},
\end{split}
\end{equation} 
where $h_{-i_b}>0$ and $h_{i_f}>0$,  or $\Delta$ is odd.  Its $\ell_1$ norm condition is proved in the literature above as
\begin{equation}\label{eq:FIR_noncausal_L1}
	\sum_{i_b=1}^{n_b}|h_{-i_b}|\rho^{-2i_b} + \sum_{i_f=1}^{n_f}|h_{i_f}|< h_0,  
\end{equation}
which in turn proves (\ref{eq:FIR_noncausal_L11}).

As the variable conversion is unique, the multipliers $M(\rho,z)$ in (\ref{eq:FIR_noncausal}) and $M(z)$ in (\ref{eq:FIR_noncausal1}) are equivalent and  belong to the same subset $\mathcal{M}_{\rho}$.  Nevertheless, the main issue of the noncausal multiplier $M(\rho,z)$ is the factorisation.

First, for causal multipliers, the factorisations (\ref{eq:zf_factorisation1}), (\ref{eq:zf_factorisation2}) and (\ref{eq:zf_factorisation3}) are all valid with $z$ replaced by $\rho z$. After that the LMI (\ref{eq:causal_LMI}) is obtained with the  KYP lemma.  In addition, when the factorisation (\ref{eq:phi_pz}) is used, the LMI is in the form (\ref{eq:iqc_lmi}).

Second, for anticausal multipliers, the factorisation (\ref{eq:zf_factorisation1}) is invalid. The other two factorisations by lifting are discussed.

The factorisation (\ref{eq:zf_factorisation2}) and  (\ref{eq:zf_factorisation3}) lead to (\ref{eq:factorisation_rho}) with anticausal multipliers by setting $\pmb{h}_{-i}=\pmb{0}$ and $\pmb{h}_{-i_b}=\pmb{0}$, respectively, and replacing $z$ by $\frac{z}{\rho}$ in $\Psi_{2,3}(z)$. The  corresponding state-space representation becomes
\begin{gather*}
\Psi_{2,3}(\rho,z) \sim  \begin{bmatrix}  \rho A_\Psi & \rho B_{\Psi_1} & \rho B_{\Psi_2}\\ C_\Psi & D_{\Psi_1} & D_{\Psi_2} \end{bmatrix}.  
\end{gather*}

Notice that $\Psi_{3}(\rho,z)$ here is only valid for anticausal multipliers, while it for noncausal multipliers is (\ref{eq:phi_pz}).

With this factorisation,  the replacements $A_{\Psi} \to \rho A_\Psi $, $B_{\Psi} \to \rho B_\Psi $ are taken in  (\ref{eq:iqc_lmi}) to obtain the LMI in the anticausal case.  

The factorisation (\ref{eq:phi_pz}) leads to  (\ref{eq:factorisation_rho}) with anticausal multipliers by setting $\pmb{h}_{-i_b}=\pmb{0}$. Similar to the causal case, the LMI (\ref{eq:iqc_lmi}) will be obtained.

Finally, for noncausal multipliers, the factorisation is more restricted, because it is inconsistent for the causal part and anticausal part. Using the factorisation (\ref{eq:zf_factorisation2}) as an example, it is impossible to replace the variable $z$ in  $\Psi_2$ to $\rho z $ for causal multipliers and to $\frac{z}{\rho}$ for anticausal multipliers simultaneously.  Therefore, the decoupling of causal and anticausal part is significant, and only the factorisation  (\ref{eq:phi_pz}) is valid. Then,  the LMI is (\ref{eq:iqc_lmi}).

In short,  by this approach, the factorisation is more restricted. The factorisation (\ref{eq:factorisation_rho}) with $\left(\Psi_3(\rho,z), K_{p,3}\right)$ is an unified structure in the $\rho$-IQC framework.

\section{Numerical results}\label{sec:results}
%\subsection{Results by different multipliers and factorisations}
In this section, we compare the results by causal, anticausal and noncausal multipliers with different forms and factorisations. The examples are listed in Table \ref{tb:example_plants} with other related information. Here, we set the causal step and anticausal step to be equal in noncausal multipliers ($n_b=n_f=n_z$).  According to the preliminary study, higher order multipliers may not lead to less-conservative results, so the step in each example was tuned in advance to reduce the conservatism.

\begin{table}[ht]
	\centering
	\begin{tabular}{ |c|c|c|c| } 
		\hline
		Ex & $G(z)$  &$K$ &$n_z$ \\ \hline
		1 & $-\frac{1}{z-0.4}$& $1$ & $1$ \\  \hline
		2 & $\frac{2z -1}{20z^2-10z+10}$& $9$ & $20$\\ \hline
		3 & $-\frac{10z^2+19z+9}{100z^3 -80z^2+17z-1}$&$3$  & $30$\\ \hline
		4 & $-\frac{0.1z}{z^2-1.8z+0.81}$& $12$  & $20$\\ \hline
	\end{tabular}
	\caption{Examples}
	\label{tb:example_plants}
\end{table}

The best estimates of $\overline{\rho^*_{\{G,K\}}}$ are obtained by the bisection search from the initial range $(max(|eig(A)|), 1)$, where the lower boundary is set to ensure the stability of the scaled plant $G(\rho z)$. The software package CVX with the solver sdpt3 \cite{cvx:2014, cvx:2008} is used to solve the LMI. The results are demonstrated in Table \ref{tb:convergence_results}, where the best result of each example is in bold. 

\begin{table}[ht]
	\centering
	\begin{tabular}{ |c|c|c|c|c|c| } 
		\hline
		 &   & \multicolumn{4}{|c|}{$\overline{\rho^*_{\{G,K\}}}$} \\ \cline{3-6}
			
		Ex &  $\underline{\rho^*_{\{G,K\}}}$ & C. (\ref{eq:FIR_noncausal}) &  AC. (\ref{eq:FIR_noncausal})  &  NC. (\ref{eq:FIR_noncausal1})  &   NC. (\ref{eq:FIR_noncausal})\\ 
		
	 	&   & with (\ref{eq:zf_factorisation2}) &  with (\ref{eq:zf_factorisation2}) & with (\ref{eq:zf_factorisation2})  & with (\ref{eq:phi_pz}) \\ \hline

		1 & $0.6$ & $0.600037$ &\pmb{$0.600024$} &\pmb{$0.600024$}&\pmb{$0.600024$}\\  \hline
		1a & $0.6$ & $0.600281$ &\pmb{$0.600024$} &\pmb{$0.600024$}&\pmb{$0.600024$}\\  \hline
		2 & $0.974679$& $0.978485$ & $0.975044$ &\pmb{$0.974758$}&$0.974794$\\ \hline
		2a & $0.974679$& $0.998474$ & $0.975044$ &\pmb{$0.974758$}&$0.974794$\\ \hline
		3 & $0.975367$& $0.976437$ & \pmb{$0.975525$}&$0.975815$ &$0.975891$\\ \hline
		3a & $0.975367$& $0.976501$ & \pmb{$0.975769$}&$0.975891$ &$0.975891$\\ \hline
		4 & $0.9$& $0.991760$ & invalid &\pmb{$0.990723$}&\pmb{$0.990723$}\\ \hline
		4a & $0.9$& $0.992794$ & invalid &\pmb{$0.992529$}&\pmb{$0.992529$}\\ \hline
	\end{tabular}
	\caption{Best estimates of $\overline{\rho^*_{\{G,K\}}}$ by causal (C.), anticausal (AC.) and noncausal (NC.) multipliers with different factorisations compared with $\underline{\rho^*_{\{G,K\}}}$}
	\label{tb:convergence_results}
\end{table}

In the table, Ex. 1-4 are with odd uncertainties while Ex. 1a-4a are with non-odd uncertainties. As demonstrated in the table, all the results are valid as they are larger than the theoretical bound $\rho^*$. For the systems that break the Kalman conjecture, such as Ex. 4(a), $\underline{\rho^*_{\{G,K\}}}$ is meaningless. Otherwise, the best results are  close to the theoretical value. 

First, as for different types of uncertainties, the optimal convergence rates with odd uncertainties are no greater than those with general uncertainties. which is a natural result and can be reflected by the more strict $\ell_1$ norm conditions on multipliers for general uncertainties. 

Then, as for different multipliers, the introduction of anticausal steps are efficient to reduce the conservatism (also see \cite[Fig. 2]{Freeman:2018}). Meanwhile, noncausal multipliers should be less conservative than causal and anticausal multipliers in general. Nevertheless, due to possible numerical problems with noncausal multipliers in MATLAB (e.g. large computational load, complex matrices), anticausal multipliers  provide the optimal results in Ex. 1(1a) and 3(3a).  However, anticausal multipliers are conservative when searching the maximum slope $K$ for  $\ell_2$-stability (also see \cite[Table \rom{2}, \rom{3}]{Shuai:2014}). For example, in Ex. 4(4a), the anticausal multiplier is not sufficient for  stability when $K=12$, while causal and noncausal multipliers are sufficient. In this case, noncausal multipliers are less conservative. 

Moreover, the results by the two structures of noncausal multipliers $M(z)$ and $M(\rho, z)$ are almost the same, as they are equivalent. Nevertheless, the computational load with $M(\rho, z)$ would be larger due to its large matrices in the factorisation. 

\begin{figure}[ht]
	\centering
	\includegraphics[width=\linewidth]{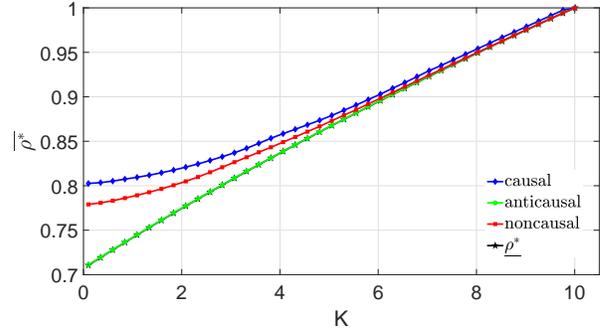}
	\caption{$\overline{\rho^*_{\{G,K\}}}$ with different $K$ and multipliers of Ex.2}
	\label{fig:p_ex2}
\end{figure}

\begin{figure}[ht]
	\centering
	\includegraphics[width=\linewidth]{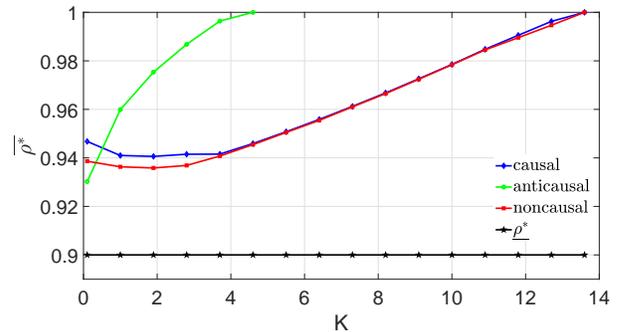}
	\caption{$\overline{\rho^*_{\{G,K\}}}$ with different $K$ and multipliers of Ex.4}
	\label{fig:p_ex4}
\end{figure}

The similar conclusions are also inflected by the convergence rates with uncertainties that have different maximum slopes in Fig. \ref{fig:p_ex2},  \ref{fig:p_ex4}.


As illustrated in the figures, anticausal multipliers achieve tight bounds of the convergence rates when they are sufficient for  exponential stability; noncausal multipliers are efficient to achieve less conservative results when the closed-loop system is close to  instability. However, causal multipliers are conservative in general.

In summary, for a general plant whose properties are unknown, it is more reliable to use noncausal multipliers to obtain less conservative bounds of the convergence rates.     Meanwhile, the choice of a specific structure of noncausal multipliers is not crucial as for the conservatism in result.  On the contrary, if the given plant verifies the Kalman conjecture and the feedback uncertainty is odd, Conjecture \ref{cj:kc_convergence} seems true, and the multiplier techniques are not necessary.

\section{Conclusion}\label{sec:conclusion}
In this technical note, we reviewed the stability concepts of Lur'e systems, where exponential stability with the convergence rate $\rho$ is linked with an extension form of  $\ell_2$-stability, defined as $\ell_2^{\rho}$ stability.   
On the other hand, we extended the literature results on causal FIR Zames-Falb multipliers to anticausal and noncausal cases. The Zames-Falb IQC and $\rho$-IQC are shown equivalent, where the multipliers can be converted by changing the variables. However, the factorisation of the Zames-Falb $\rho$-IQC is restricted, especially with noncausal multipliers. Hence, an unified factorisation is provided for both Zames-Falb IQC and $\rho$-IQC.
Furthermore, the numerical examples indicated noncausal multipliers are efficient to achieve less conservative estimation of the upper bound convergence rates in general case. However, when the system verifies the Kalman conjecture and the feedback uncertainty is odd, it is reasonable to use the theoretical bound from the linear analysis directly without any multiplier technique.

\bibliographystyle{IEEEtran}
\bibliography{refer}

\appendix \label{Sec:appendix}
\subsection{Proof of Theorem \ref{th:stability_relations} (Outline)}\label{appendix_proof}
The proof follows \cite[Theorem 6.3.46]{Vidyasagar:2002} and the theorems therein, which link  internal  stability to  input-output $\ell$-stability in continuous time. Here we relate the relations with the convergence rate in discrete time. The argument in continuous time can be converted to discrete time trivially, which is not stressed here.	

In order to keep consistent with the notations in \cite{Vidyasagar:2002},  we consider the system $-G \sim $$\begin{bmatrix}  A\;  B;\; -C\; 0 \end{bmatrix}$. The minus signs indicate negative feedback structure in \cite[Fig. 6.4]{Vidyasagar:2002}, and the $D$ matrix can be set to $0$ without loss of generality. Moreover, because $G$ is linear and stable, the disturbance signal $f$ in Fig. \ref{fig:lure} is assumed to be zero without loss of generality.  Then, the feedback interconnection in Fig. \ref{fig:lure} can be expressed as
\begin{equation}\label{eq:system_proof}
x_{k+1}=Ax_k+Bg_k+B\Delta(v_k),\quad y_k=-Cx_k.
\end{equation}

\textit{Necessity}:
The expression (\ref{eq:system_proof}) satisfies the general form (6.3.7) in \cite{Vidyasagar:2002}. In addition, the conditions (6.3.16, 17) are satisfied, because the state-space matrices of $G$ are finite  and $\Delta$ is slope-restricted. 

Following the proof and notations in \cite[Theorem 6.3.15]{Vidyasagar:2002}, (6.3.19) implies the exponentially stability condition (\ref{eq:exponential_stable}) of the unforced system where $c=\frac{\beta}{\alpha}$, $\rho=\frac{1}{2\beta^2}$ ($0<\alpha<\beta$ are constant). Then, in (6.3.32), $\|x_k\|\le \frac{W_k}{\alpha}$, where $W_k\le h_k$, and $h_k$ is the output of a first order system with input $\|g_k\|$ and pole $\frac{-1}{2\beta^2}$. In other words, the exponential rate limited by this transfer function is $\rho=\frac{1}{2\beta^2}$. Henceforth, let the input to this transfer function is $g\in \ell_2^{\rho}$, then the solution $h\in \ell_2^{\rho}$, thus $W\in \ell_2^{\rho}$  and $x\in \ell_2^{\rho}$. Moreover, according to (6.3.33), $v\equiv y\in \ell_2^{\rho}$. Next, because $\Delta$ is memoryless  and slope-restricted, and $g\in \ell_2^{\rho}$, the signal $w\in \ell_2^{\rho}$. This proves that the forced system is  $\ell_2^{\rho}$-stable.

%\textit{If}: In the proof of \cite[Theorem 6.3.46]{Vidyasagar:2002}, the closed-loop system (\ref{eq:system_proof}) is reachable and uniformly observable when $G$ is controllable and observable. Then, according to the proof in \cite[Theorem 6.3.39]{Vidyasagar:2002}, there exists an input $g_k$ with $g_k=0$ for $k>k^*$ (thus $g\in\ell_2^{\rho}$) to achieve $x_{k^*}$ with zero initial state. Therefore, the state of the unforced system $x_{u,k}$ with the initial state $x_{k^*}$ and zero input is equal to the state of the forced system $x_{f,k^*+k}$ with the input $g$ and zero initial state. Hereafter, because  $g\in\ell_2^{\rho}$ and $y\in\ell_2^{\rho}$, (6.3.41) can be rewritten as
%\begin{equation}\label{eq:y1_convergence}
%\sum_{k}^{\infty}\|y_k \rho^{-k}\|^2\to 0, \quad \text{as} \; k\to\infty.
%\end{equation} 
%
%In addition, (6.3.34) implies 
%\begin{equation}\label{eq:y2_convergence}
%\sum_{k}^{\infty}\|y_k \rho^{-k}\|^2\ge \alpha \|x_k \rho^{-k}\|, \quad \forall k\ge k^*.
%\end{equation} 
%
%The relations (\ref{eq:y1_convergence}) and (\ref{eq:y2_convergence}) show that  $x_k \rho^{-k}\to 0$ as $k\to \infty$. Finally, by Remark \ref{rm:exponential_convergence}, the unforced system is exponentially stable with rate $\rho$.

\subsection{State-space representation of $\Psi_2(z)$ and $\Psi_3(\rho,z)$} \label{appendix_ss}

The state-space representation of $\Psi_2(z)$ in (\ref{eq:zf_factorisation2}) is given below. 
{\begin{equation*}
	\begin{split}
	A_{2\Psi}=\begin{bmatrix}
	A_S & \pmb{0}\\ \pmb{0} & A_S
	\end{bmatrix}_{(2n\times 2n)},&\quad
	B_{2\Psi}=\begin{bmatrix}
	B_S & \pmb{0}\\ \pmb{0} & B_S
	\end{bmatrix}_{(2n \times 2)}, \\
	C_{2\Psi}=\begin{bmatrix}
	C_{S} & \pmb{0}\\ \pmb{0} & C_{S}
	\end{bmatrix}_{((2n+2) \times 2n)},& \quad
	D_{2\Psi}=\begin{bmatrix}
	D_{S} & \pmb{0}\\ \pmb{0} & D_{S}
	\end{bmatrix}_{((2n+2) \times 2)}, 
	\end{split}
	\end{equation*}
	where $A_S$, $B_S$, $C_{S}$  and $D_{S}$ are 
	\begin{equation*}
	\begin{split}
	A_S=\begin{bmatrix}
	\pmb{0} & 0\\ \pmb{I}_{(n-1)} & \pmb{0}
	\end{bmatrix}_{(n \times n)},&\quad 
	B_S=\begin{bmatrix}
	1 \\ \pmb{0}_{(n-1) \times 1}
	\end{bmatrix}_{(n \times 1)}, \\
	C_{S}=\begin{bmatrix}
	\pmb{0} \\ \pmb{I}_{(n)}
	\end{bmatrix}_{((n+1) \times n)},&\quad 
	D_{S}=\begin{bmatrix}
	1 \\ \pmb{0}_{n \times 1}
	\end{bmatrix}_{((n+1) \times 1)}.
	\end{split}
	\end{equation*}

Then, we consider $\Psi_3(\rho,z)$ in (\ref{eq:phi_pz}). Let $n_b \ne n_f$, and $\hat{n}=\max\{n_b,nf\}$. The state-space representation is given below. 
{\begin{equation*}
\begin{split}
A_{3\Psi}=\frac{1}{\rho}A_{2\Psi (2\hat{n} \times 2\hat{n})},&\quad
B_{3\Psi}=B_{2\Psi (2\hat{n} \times 2)}, \\
C_{3\Psi}=\begin{bmatrix}
C_{S1} & \pmb{0}\\ \pmb{0} & C_{S1}\\C_{S2} & \pmb{0}\\ \pmb{0} & C_{S2}
\end{bmatrix},& \;
D_{3\Psi}=\begin{bmatrix}
D_{S ((n_b+1) \times 1)} & \pmb{0}\\ \pmb{0} & D_{S ((n_b+1) \times 1)} \\D_{S ((n_f+1) \times 1)}  & \pmb{0}\\ \pmb{0} & D_{S ((n_f+1) \times 1)} 
\end{bmatrix}.
\end{split}
\end{equation*}

The matrices $C_{S1}$  and $C_{S2}$ are 

{\small{\begin{equation*}
	C_{S1}=\arraycolsep1pt \begin{bmatrix}
	\pmb{0} & \pmb{0}^{\square}\\ \frac{1}{\rho}\pmb{I}_{(n_b)} & \pmb{0}^{\square}
	\end{bmatrix}_{((n_b+1)\times\hat{n})},
	C_{S2}=\begin{bmatrix}
	\pmb{0} & \pmb{0}^{\triangle}\\ diag\left(\rho^{2i_f-1}\right) & \pmb{0}^{\triangle}
	\end{bmatrix}_{((n_f+1)\times\hat{n})}.
	\end{equation*}}}

The matrices $\pmb{0}^{\square}$ are removed when $n_b>n_f$; the matrices $\pmb{0}^{\triangle}$ are removed when $n_b<n_f$; both of the matrices $\pmb{0}^{\square}$ and $\pmb{0}^{\triangle}$ are removed when $n_b=n_f$.

\end{document}